\documentclass[12pt]{article}
\newcommand{\be}{\begin{equation}}
\newcommand{\ee}{\end{equation}}
\newcommand{\bea}{\begin{eqnarray}}
\newcommand{\eea}{\end{eqnarray}}
\newcommand{\ba}{\begin{array}}
\newcommand{\ea}{\end{array}}

\newcommand{\bc}{\begin{center}}
\newcommand{\ec}{\end{center}}
\newcommand{\ben}{\begin{enumerate}}
\newcommand{\een}{\end{enumerate}}
\newcommand{\bfi}{\begin{figure}}
\newcommand{\efi}{\end{figure}}

\newcommand{\bq}{\begin{quote}}
\newcommand{\eq}{\end{quote}}
\newcommand{\bqu}{\begin{quotation}}
\newcommand{\equ}{\end{quotation}}
\newenvironment{emphit}{\begin{itemize}}{\end{itemize}}
\newcommand{\bemp}{\begin{emphit}}
\newcommand{\eemp}{\end{emphit}}

\newcommand{\bt}{\begin{tabular}}
\newcommand{\et}{\end{tabular}}

\newtheorem{myth}{Theorem}[section]
\newtheorem{mylem}{Lemma}[section]
\newtheorem{mycor}{Corollary}[section]

\newtheorem{mydef}{Definition}[section]
\newtheorem{myrem}{Remark}[section]
\newtheorem{myexam}{Example}[section]

\begin{document}
\date{}
\title{Weak Inverse
Property Loops and Some Isotopy-Isomorphy Properties \footnote{2000
Mathematics Subject Classification. Primary 20NO5 ; Secondary 08A05}
\thanks{{\bf Keywords and Phrases :}weak inverse property, cross inverse property, isotopism}}
\author{T. G. Jaiy\'e\d ol\'a\thanks{All correspondence to be addressed to this author.} \\
Department of Mathematics,\\
Obafemi Awolowo University,\\
Ile Ife 220005, Nigeria.\\
jaiyeolatemitope@yahoo.com\\tjayeola@oauife.edu.ng \and
J. O. Ad\'en\'iran \\
Department of Mathematics,\\
University of Abeokuta, \\
Abeokuta 110101, Nigeria.\\
ekenedilichineke@yahoo.com\\
adeniranoj@unaab.edu.ng} \maketitle
\begin{abstract}
Two distinct isotopy-isomorphy conditions, different from those of
J. M. Osborn and Wilson's condition, for a weak inverse property
loop(WIPL) are shown. Only one of them characterizes
isotopy-isomorphy in WIPLs while the other is just a sufficient
condition for isotopy-isomorphy. Under the sufficient condition
called the ${\cal T}$ condition, Artzy's result that isotopic cross
inverse property loops(CIPL) are isomorphic is proved for WIP loops.
\end{abstract}

\section{Introduction}
\paragraph{}
Michael K. Kinyon \cite{phd33} gave a talk on Osborn Loops and
proposed the open problem : "Is every Osborn Loop universal?" which
is obviously true for universal WIP loops and universal CIP loops. A
popular isotopy-isomorphy condition in loops is the Wilson's
identity(\cite{phd91}) and a loop obeying it is called a Wilson's
loop by Goodaire and Robinson \cite{phd48, phd91} and they proved
that a loop is a Wilson loop if and only if it is a conjugacy closed
loop(CC-loop) and a WIPL. Our aim in this work is to prove some
isotopy-isomorphy conditions, different from those of J. M. Osborn
\cite{phd89} and Wilson \cite{phd90}(i.e the loop not necessarily a
CC-loop), for a WIPL and see if the result of Artzy \cite{phd30}
that isotopic CIP loops are isomorphic is true for WIP loops or some
specially related WIPLs(i.e special isotopes). But before these, we
shall take few basic definitions and concepts in loop theory which
are needed here.

Let $L$ be a non-empty set. Define a binary operation ($\cdot $) on
$L$ : If $x\cdot y\in L$ for all $x, y\in L$, $(L, \cdot )$ is
called a groupoid. If the system of equations ; $a\cdot x=b$ and
$y\cdot a=b$ have unique solutions for $x$ and $y$ respectively,
then $(L, \cdot )$ is called a quasigroup. Furthermore, if there
exists a unique element $e\in L$ called the identity element such
that for all $x\in L$, $x\cdot e=e\cdot x=x$, $(L, \cdot )$ is
called a loop. For each $x\in L$, the elements $x^\rho ,x^\lambda\in
L$ such that $xx^\rho=e=x^\lambda x$ are called the right, left
inverses of $x$ respectively. $L$ is called a weak inverse property
loop (WIPL) if and only if it obeys the weak inverse property (WIP);
$xy\cdot z=e$ implies $x\cdot yz=e$ for all $x,y,z\in L$ while $L$
is called a cross inverse property loop (CIPL) if and only if it
obeys the cross inverse property (CIP); $xy\cdot x^\rho=y$.

According to \cite{phd31}, the WIP is a generalization of the CIP.
The latter was introduced and studied by R. Artzy \cite{phd45} and
\cite{phd30} while the former was introduced by J. M. Osborn
\cite{phd89} who also investigated the isotopy invariance of the
WIP. Huthnance Jr. \cite{phd44} did so as well and proved that the
holomorph of a WIPL is a WIPL.  A loop property is called
universal(or at times a loop is said to be universal relative to a
particular property) if the loop has the property and every loop
isotope of such a loop possesses such a property. A universal WIPL
is called an Osborn loop in Huthnance Jr. \cite{phd44} but this is
different from the Osborn loop of Kinyon \cite{phd33} and Basarab.
The Osborn loops of Kinyon and Basarab were named generalised
Moufang loops or M-loops by Huthnance Jr. \cite{phd44} where he
investigated the structure of their holomorphs while Basarab
\cite{phd46} studied Osborn loops that are G-loops. Also,
generalised Moufang loops or M-loops of Huthnance Jr. are different
from those of Basarab. After Osborn's study of universal WIP loops,
Huthnance Jr. still considered them in his thesis and did an
elaborate study by comparing the similarities between properties of
Osborn loops(universal WIPL) and generalised Moufang loops. He was
able to draw conclusions that the latter class of loops is large
than the former class while in a WIPL the two are the same.

But in this present work, two distinct isotopy-isomorphy conditions,
different from that of Osborn \cite{phd89} and Wilson \cite{phd90},
for a weak inverse property loop(WIPL) are shown. Only one of them
characterizes isotopy-isomorphy in WIPLs while the other is just a
sufficient condition for isotopy-isomorphy. Under the sufficient
condition called the ${\cal T}$ condition, Artzy \cite{phd30} result
that isotopic cross inverse property loops(CIPL) are isomorphic is
proved for WIP loops.
\section{Preliminaries}
\begin{mydef}
Let $(L, \cdot )$ and $(G, \circ )$ be two distinct loops. The
triple $\alpha=(U, V, W) : (L, \cdot )\to (G, \circ )$ such that $U,
V, W : L\to G$ are bijections is called a loop isotopism
$\Leftrightarrow~ xU\circ yV=(x\cdot y)W~\forall ~x, y\in L$. Hence,
$L$ and $G$ are said to be isotopic whence, $G$ is an isotope of
$L$.
\end{mydef}
\begin{mydef}
Let $L$ be a loop. A mapping $\alpha\in S(L)$(where $S(L)$ is the
group of all bijections on $L$) which obeys the identity
$x^\rho=[(x\alpha)^\rho ]\alpha$ is called a weak right inverse
permutation. Their set is represented by $S_\rho (L)$.

Similarly, if $\alpha$ obeys the identity
$x^\lambda=[(x\alpha)^\lambda
 ]\alpha$ it is called a weak left inverse permutation. Their set is represented by $S_\rho (L)$

If $\alpha $ satisfies both, it is called a weak inverse
permutation. Their set is represented by $S'(L)$.

It can be shown that $\alpha\in S(L)$ is a weak right inverse if and
only if it is a weak left inverse permutation. So, $S'(L)=S_\rho
(L)=S_\lambda (L)$.
\end{mydef}

\begin{myrem}\label{0:1}
Every permutation of order 2 that preserves the right(left) inverse
of each element in a loop is a weak right (left) inverse
permutation.
\end{myrem}

\begin{myexam}\label{0:3}
If $L$ is an extra loop, the left and right inner mappings $L(x,y)$
and $R(x,y)~\forall~x,y\in L$ are automorphisms of orders 2
(\cite{phd36}). Hence, they are weak inverse permutations by
Remark~\ref{0:1}
\end{myexam}
\paragraph{}Throughout, we shall employ the use of the bijections;
$J_\rho~:~x\mapsto x^\rho$, $J_\lambda~:~x\mapsto x^\lambda$,
$L_x~:~y\mapsto xy$ and $R_x~:~y\mapsto yx$ for a loop and the
bijections; $J_\rho'~:~x\mapsto x^{\rho'}$, $J_\lambda'~:~x\mapsto
x^{\lambda'}$, $L_x'~:~y\mapsto xy$ and $R_x'~:~y\mapsto yx$ for its
loop isotope. If the identity element of a loop is $e$ then that of
the isotope shall be denoted by $e'$.

\begin{mylem}\label{0:2}
In a loop, the set of weak inverse permutations that commute form an
abelian group.
\end{mylem}

\begin{myrem}\label{0:4}
Applying Lemma~\ref{0:2} to extra loops and considering
Example~\ref{0:3}, it will be observed that in an extra loop $L$,
the Boolean groups $Inn_\lambda(L),Inn_\rho\le S'(L)$ .
$Inn_\lambda(L)$ and $Inn_\rho(L)$ are the left and right inner
mapping groups respectively. They have been investigated in
\cite{phd35} and \cite{phd36}. This deductions can't be drawn for
CC-loops despite the fact that the left (right) inner mappings
commute and are automorphisms. And this is as a result of the fact
that the left(right) inner mappings are not of exponent 2.
\end{myrem}
\begin{mydef}(${\cal T}$-condition)

Let $(G,\cdot )$ and $(H,\circ )$ be two distinct loops that are
isotopic under the triple $(A,B,C)$. $(G,\cdot )$ obeys the ${\cal
T}_1$ condition if and only if $A=B$. $(G,\cdot )$ obeys the ${\cal
T}_2$ condition if and only if $J_\rho'=C^{-1}J_\rho B=A^{-1}J_\rho
C$. $(G,\cdot )$ obeys the ${\cal T}_3$ condition if and only if
$J_\lambda'=C^{-1}J_\lambda A=B^{-1}J_\lambda C$. So, $(G,\cdot )$
obeys the ${\cal T}$ condition if and only if it obey ${\cal T}_1$
and ${\cal T}_2$ conditions or ${\cal T}_1$ and ${\cal T}_3$
conditions since ${\cal T}_2\equiv {\cal T}_3$.
\end{mydef}
It must here by be noted that the ${\cal T}$-conditions refer to a
pair of isotopic loops at a time. This statement might be omitted at
times. That is whenever we say a loop $(G,\cdot )$ has the ${\cal
T}$-condition, then this is relative to some isotope $(H,\circ )$ of
$(G,\cdot )$

\begin{mylem}\label{1:1}(\cite{phd3})
Let $L$ be a loop. The following are equivalent.
\begin{enumerate}
\item $L$ is a WIPL
\item $y(xy)^\rho =x^\rho~\forall~x,y\in L$.
\item $(xy)^\lambda x=y^\lambda~\forall~x,y\in L$.
\end{enumerate}
\end{mylem}
\begin{mylem}\label{1:2}
Let $L$ be a loop. The following are equivalent.
\begin{enumerate}
\item $L$ is a WIPL
\item $R_yJ_\rho L_y=J_\rho~\forall~y\in L$.
\item $L_xJ_\lambda R_x=J_\lambda~\forall~x\in L$.
\end{enumerate}
\end{mylem}
\section{Main Results}
\begin{myth}\label{1:3}
Let $(G,\cdot )$ and $(H,\circ )$ be two distinct loops that are
isotopic under the triple $(A,B,C)$.
\begin{enumerate}
\item If the pair of $(G,\cdot )$ and $(H,\circ )$ obey the ${\cal T}$
condition, then $(G,\cdot )$ is a WIPL if and only if $(H,\circ )$
is a WIPL.
\item If  $(G,\cdot )$ and $(H,\circ )$ are WIPLs, then $J_\lambda R_xJ_\rho B=CJ_\lambda' R_{xA}'J_\rho'$ and $J_\rho
L_xJ_\lambda A=CJ_\rho' L_{xB}'J_\lambda'$ for all $x\in G$.
\end{enumerate}
\end{myth}
{\bf Proof}\\
\begin{enumerate}
\item $(A, B,C)~:~G\rightarrow H$ is an isotopism $\Leftrightarrow xA\circ
yB=(x\cdot y)C\Leftrightarrow yBL_{xA}'=yL_xC\Leftrightarrow
BL_{xA}'=L_xC\Leftrightarrow L_{xA}'=B^{-1}L_xC\Leftrightarrow$
\begin{equation}\label{eq:1}
L_x=BL_{xA}'C^{-1}
\end{equation}
Also, $(A, B,C)~:~G\rightarrow H$ is an isotopism $\Leftrightarrow
xAR_{yB}'=xR_yC\Leftrightarrow AR_{yB}'=R_yC\Leftrightarrow
R_{yB}'=A^{-1}R_yC\Leftrightarrow$
\begin{equation}\label{eq:2}
R_y=AR_{yB}'C^{-1}
\end{equation}
Applying (\ref{eq:1}) and (\ref{eq:2}) to Lemma~\ref{1:2}
separately, we have : $R_yJ_\rho L_y=J_\rho $, $L_xJ_\lambda
R_x=J_\lambda \Rightarrow (AR_{xB}'C^{-1})J_\rho
(BL_{xA}'C^{-1})=J_\rho$, $(BL_{xA}'C^{-1})J_\lambda
(AR_{xB}'C^{-1})=J_\lambda\Leftrightarrow AR_{xB}'(C^{-1}J_\rho
B)L_{xA}'C^{-1}=J_\rho$, $BL_{xA}'(C^{-1}J_\lambda
A)R_{xB}'C^{-1}=J_\lambda\Leftrightarrow$
\begin{equation}\label{eq:3}
R_{xB}'(C^{-1}J_\rho B)L_{xA}'=A^{-1}J_\rho
C,~L_{xA}'(C^{-1}J_\lambda A)R_{xB}'=B^{-1}J_\lambda C.
\end{equation}

Let $J_\rho'=C^{-1}J_\rho B=A^{-1}J_\rho C$,
$J_\lambda'=C^{-1}J_\lambda A=B^{-1}J_\lambda C$. Then, from
(\ref{eq:3}) and by Lemma~\ref{1:2}, $H$ is a WIPL if $xB=xA$ and
$J_\rho'=C^{-1}J_\rho B=A^{-1}J_\rho C$ or $xA=xB$ and
$J_\lambda'=C^{-1}J_\lambda A=B^{-1}J_\lambda C \Leftrightarrow B=A$
and $J_\rho'=C^{-1}J_\rho B=A^{-1}J_\rho C$ or $A=B$ and
$J_\lambda'=C^{-1}J_\lambda A=B^{-1}J_\lambda C \Leftrightarrow A=B$
and $J_\rho'=C^{-1}J_\rho B=A^{-1}J_\rho C$ or
$J_\lambda'=C^{-1}J_\lambda A=B^{-1}J_\lambda C$. This completes the
proof of the forward part. To prove the converse, carry out the same
procedure, assuming the ${\cal T}$ condition and the fact that
$(H,\circ )$ is a WIPL.
\item If $(H,\circ )$ is a WIPL, then
\begin{equation}\label{eq:3.1}
R_y'J_\rho' L_y'=J_\rho',~\forall~y\in H
\end{equation}
while since $G$ is a WIPL,
\begin{equation}\label{eq:3.2}
R_xJ_\rho L_x=J_\rho~\forall~x\in G.
\end{equation}
The fact that $G$ and $H$ are isotopic implies that
\begin{equation}\label{eq:3.3}
L_x=BL_{xA}'C^{-1}~\forall~x\in G~and
\end{equation}
\begin{equation}\label{eq:3.4}
R_x=AR_{xB}'C^{-1}~\forall~x\in G.
\end{equation}
From (\ref{eq:3.1}),
\begin{equation}\label{eq:3.5}
R_y'=J_\rho' L_y'^{-1}J_\lambda'~\forall~y\in H~and
\end{equation}
\begin{equation}\label{eq:3.6}
L_y'=J_\lambda' R_y'^{-1}J_\rho'~\forall~y\in H
\end{equation}
while from (\ref{eq:3.2}),
\begin{equation}\label{eq:3.7}
R_x=J_\rho L_x^{-1}J_\lambda~\forall~x\in G~and
\end{equation}
\begin{equation}\label{eq:3.8}
L_x=J_\lambda R_x^{-1}J_\rho~\forall~x\in G.
\end{equation}
So, using (\ref{eq:3.6}) and (\ref{eq:3.8}) in (\ref{eq:3.3}) we get
\begin{equation}\label{eq:3.9}
J_\lambda R_xJ_\rho B=CJ_\lambda' R_{xA}'J_\rho'~\forall~x\in G
\end{equation}
while using (\ref{eq:3.5}) and (\ref{eq:3.7}) in (\ref{eq:3.4}) we
get
\begin{equation}\label{eq:3.10}
J_\rho L_xJ_\lambda A=CJ_\rho' L_{xB}'J_\lambda'~\forall~x\in G.
\end{equation}
\end{enumerate}
\begin{myrem}
In Theorem~\ref{1:3}, a loop is a universal WIPL under the ${\cal
T}$ condition. But the converse of this is not true. This can be
deduced from a counter example.
\paragraph{Counter Example}
Let $G=\{0,1,2,3,4,\}$. From the Table~\ref{wipl1}, $(G,\cdot )$ is
a WIPL.
\begin{center}
\begin{table}[tbp]
\begin{tabular}{|c||c|c|c|c|c|c|}
\hline
$\cdot $ & 0 & 1 & 2 & 3 & 4 \\
\hline  \hline
0 & 0 & 1 & 2 & 3 & 4 \\
\hline
1 & 1 & 3 & 0 & 4 & 2 \\
\hline
2 & 2 & 0 & 4 & 1 & 3 \\
\hline
3 & 3 & 4 & 1 & 2 & 0 \\
\hline
4 & 4 & 2 & 3 & 0 & 1 \\
\hline
\end{tabular}
\caption{A commutative weak inverse property loop}\label{wipl1}
\end{table}
\end{center}
Let \begin{displaymath} A= \left( \begin{array}{ccccc}
0 & 1 & 2 & 3 & 4\\
1 & 3 & 0 & 4 & 2
\end{array} \right )~
\textrm{and} ~ B= \left( \begin{array}{ccccc}
0 & 1 & 2 & 3 & 4\\
2 & 0 & 4 & 1 & 3
\end{array} \right )
\end{displaymath}
Then, $(A,B,I)$ is an isotopism from $(G,\cdot )$ to itself. But
$A\ne B$ so the ${\cal T}$ condition does not hold for $(G,\cdot )$.
\end{myrem}

\begin{myth}\label{1:3.2}
Let $(G,\cdot )$ be a WIPL with identity element $e$ and $(H,\circ
)$ be an arbitrary loop isotope of $(G,\cdot )$ with identity
element $e'$ under the triple $\alpha =(A,B,C)$. If $(H,\circ )$ is
a WIPL then
\begin{enumerate}
\item $(G,\cdot )\cong^{{}^C} (H,\circ )\Leftrightarrow (J_\rho
L_bJ_\lambda , J_\lambda R_aJ_\rho ,I)\in AUT(G,\cdot )$ where
$a=e'A^{-1}, b=e'B^{-1}$. Hence, $(J_\lambda R_aJ_\rho  , J_\rho
L_bJ_\lambda,R_aL_b)\in AUT(G,\cdot )$. Furthermore, if $(G,\cdot )$
is a loop of exponent 2 then, $(R_a,L_b,R_aL_b)\in AUT(G,\cdot )$.
\item $(G,\cdot )\cong^{{}^C} (H,\circ )\Leftrightarrow (J_\rho'
L_{b'}'J_\lambda' , J_\lambda' R_{a'}'J_\rho' ,I)\in AUT(H,\circ )$
where $a'=eA, b'=eB$. Hence, $(J_\lambda' R_{a'}'J_\rho'  , J_\rho'
L_{b'}'J_\lambda',R_{a'}'L_{b'}')\in AUT(H,\circ )$. Furthermore, if
$(H,\circ )$ is a loop of exponent 2 then,
$(R_{a'}',L_{b'}',R_{a'}'L_{b'}')\in AUT(H,\cdot )$.
\item $(G,\cdot )\cong^{{}^C} (H,\circ )\Leftrightarrow (L_b,R_a,I)\in AUT(G,\cdot )$, $a=e'A^{-1},
b=e'B^{-1}$ provided $(x\cdot y)^\rho =x^\rho\cdot y^\lambda$ or
$(x\cdot y)^\lambda =x^\lambda\cdot y^\rho~\forall~x,y\in G$. Hence,
$(G,\cdot )$ and $(H,\circ )$ are isomorphic CIP loops while
$R_aL_b=I$, $ba=e$.
\item $(G,\cdot )\cong^{{}^C} (H,\circ )\Leftrightarrow (L_{b'}',R_{a'}',I)\in AUT(H,\circ )$, $a'=eA,
b'=eB$ provided $(x\circ y)^{\rho'} =x^{\rho'}\circ y^{\lambda'}$ or
$(x\circ y)^{\lambda'} =x^{\lambda'}\circ y^{\rho'}~\forall~x,y\in
H$. Hence, $(G,\cdot )$ and $(H,\circ )$ are isomorphic CIP loops
while $R_{a'}'L_{b'}'=I$, $b'a'=e'$.
\end{enumerate}
\end{myth}
{\bf Proof}\\
Consider the second part of Theorem~\ref{1:3}.
\begin{enumerate}
\item Let $y=xA$ in (\ref{eq:3.9}) and replace $y$ by $e'$. Then $J_\lambda R_{e'A^{-1}}J_\rho
B=C\Rightarrow C=J_\lambda R_aJ_\rho B\Rightarrow B=J_\lambda
R_a^{-1}J_\rho C$. Let $y=xB$ in (\ref{eq:3.10}) and replace $y$ by
$e'$. Then $J_\rho L_{e'B^{-1}}J_\lambda A=C\Rightarrow C=J_\rho
L_bJ_\lambda A\Rightarrow A=J_\rho L_b^{-1}J_\lambda C$. So, $\alpha
=(A,B,C)=(J_\rho L_b^{-1}J_\lambda C,J_\lambda R_a^{-1}J_\rho
C,C)=(J_\rho L_b^{-1}J_\lambda ,J_\lambda R_a^{-1}J_\rho
,I)(C,C,C)$. Thus, $(J_\rho L_bJ_\lambda , J_\lambda R_aJ_\rho
,I)\in AUT(G,\cdot )\Leftrightarrow (G,\cdot )\cong^{{}^C} (H,\circ
)$.

Using the results on autotopisms of WIP loops in
[Lemma~1,\cite{phd89}], $(J_\lambda R_aJ_\rho ,I,L_b),(I,J_\rho
L_bJ_\lambda,R_a)\Rightarrow (J_\lambda R_aJ_\rho  , J_\rho
L_bJ_\lambda,R_aL_b)  \in AUT(G,\cdot )$. The further conclusion
follows by breaking this
\item This is similar to (1.) above but we only need to replace $x$ by $e$
in (\ref{eq:3.9}) and (\ref{eq:3.10}).
\item This is achieved by simply breaking the autotopism in (1.) and
using the fact that a WIPL with the A. I. P. is a CIPL.
\item Do what was done in (3.) to (2.).
\end{enumerate}
\begin{mycor}\label{1:5}
Let $(G,\cdot )$ and $(H,\circ )$ be two distinct loops that are
isotopic under the triple $(A,B,C)$. If $G$ is a WIPL with the
${\cal T}$ condition, then $H$ is a WIPL :
\begin{enumerate}
\item there exists $\alpha ,\beta\in S'(G)$ i.e  $\alpha$ and $\beta $ are weak inverse
permutations and
\item $J_\rho'=J_\lambda'\Leftrightarrow J_\rho =J_\lambda$.
\end{enumerate}
\end{mycor}
{\bf Proof}\\
By Theorem~\ref{1:3}, $A=B$ and $J_\rho'=C^{-1}J_\rho B=A^{-1}J_\rho
C$ or $J_\lambda'=C^{-1}J_\lambda A=B^{-1}J_\lambda C$.
\begin{enumerate}
\item $C^{-1}J_\rho B=A^{-1}J_\rho C\Leftrightarrow J_\rho
B=CA^{-1}J_\rho C\Leftrightarrow J_\rho=CA^{-1}J_\rho CB^{-1}=
CA^{-1}J_\rho CA^{-1}=\alpha J_\rho \alpha $ where $\alpha
=CA^{-1}\in S(G,\cdot )$.
\item $C^{-1}J_\lambda A=B^{-1}J_\lambda C\Leftrightarrow J_\lambda
A=CB^{-1}J_\lambda C\Leftrightarrow J_\lambda = CB^{-1}J_\lambda
CA^{-1}=CB^{-1}J_\lambda CB^{-1}=\beta J_\lambda \beta$ where $\beta
=CB^{-1}\in S(G,\cdot )$.
\item $J_\rho'=C^{-1}J_\rho B$, $J_\lambda'=C^{-1}J_\lambda A$.
$J_\rho'=J_\lambda'\Leftrightarrow C^{-1}J_\rho B=C^{-1}J_\lambda
A=C^{-1}J_\lambda B\Leftrightarrow J_\lambda =J_\rho $.
\end{enumerate}
\begin{mylem}\label{1:6}
Let $(G,\cdot )$ be a WIPL with the ${\cal T}$ condition and
isotopic to another loop $(H,\circ )$. $(H,\circ)$ is a WIPL and $G$
has a weak inverse permutation.
\end{mylem}
{\bf Proof}\\
From the proof of Corollary~\ref{1:5}, $\alpha =\beta$, hence the
conclusion.
\begin{myth}\label{1:7}
With the ${\cal T}$ condition, isotopic WIP loops are isomorphic.
\end{myth}
{\bf Proof}\\
From Lemma~\ref{1:6}, $\alpha =I$ is a weak inverse permutation. In
the proof of Corollary~\ref{1:5}, $\alpha =CA^{-1}=I\Rightarrow
A=C$. Already, $A=B$, hence $(G,\cdot )\cong (H,\circ )$.
\begin{myrem}\label{1:3.3}
Theorem~\ref{1:3.2} and Theorem~\ref{1:7} describes isotopic WIP
loops that are isomorphic by
\begin{enumerate}
\item an autotopism in either the domain loop or the co-domain loop
and
\item the ${\cal T}$ condition(for a special case).
\end{enumerate}
These two conditions are completely different from that shown in
[Lemma~2,\cite{phd89}] and [Theorem~4,\cite{phd90}]. Furthermore, it
can be concluded from Theorem~\ref{1:3.2} that isotopic CIP loops
are not the only isotopic WIP loops that are isomorphic as earlier
shown [Theorem~1, \cite{phd30}]. In fact, isotopic CIP loops need
not satisfy Theorem~\ref{1:7}(i.e the ${\cal T}$ condition) to be
isomorphic.
\end{myrem}
\section{Conclusion and Future Study}
Karklin$\ddot{\textrm{u}}$sh and Karkli$\check{\textrm{n}}$
\cite{phd175} introduced $m$-inverse loops i.e loops that obey any
of the equivalent conditions
\begin{displaymath}
(xy)J_\rho^m\cdot xJ_\rho^{m+1}=yJ_\rho^m\qquad\textrm{and}\qquad
xJ_\lambda^{m+1}\cdot (yx)J_\lambda^m=yJ_\lambda^m.
\end{displaymath}
They are generalizations of WIPLs and CIPLs, which corresponds to
$m=-1$ and $m=0$ respectively. After the study of $m$-loops by
Keedwell and Shcherbacov \cite{phd176}, they have also generalized
them to quasigroups called $(r,s,t)$-inverse quasigroups in
\cite{phd177} and \cite{phd178}. It will be interesting to study the
universality of $m$-inverse loops and $(r,s,t)$-inverse quasigroups.
These will generalize the works of J. M. Osborn and R. Artzy on
universal WIPLs and CIPLs respectively.


\begin{thebibliography}{99}
\bibitem{phd45} R. Artzy (1955), {\it On loops with special
property}, Proc. Amer. Math. Soc. 6, 448--453.
\bibitem{phd30} R. Artzy (1959), {\it Crossed inverse and related
loops}, Trans. Amer. Math. Soc. 91, 3, 480--492.
\bibitem{phd31} R. Artzy (1960), {\it Relations between loops
identities}, Proc. Amer. Math. Soc. 11,6, 847--851.
\bibitem{phd46} A. S. Basarab (1994), {\it Osborn's G-loop}, Quasigroups and Related
Systems 1, 51--56.
\bibitem{phd91} E. G. Goodaire and D. A. Robinson (1982), {\it A
class of loops which are isomorphic to all loop isotopes}, Can. J.
Math. 34, 662--672.
\bibitem{phd48} E. G. Goodaire and D. A. Robinson (1990), {\it Some special conjugacy closed loops}, Canad. Math. Bull. 33, 73--78.
\bibitem{phd44} E. D. Huthnance Jr.(1968), {\it A theory of
generalised Moufang loops}, Ph.D. thesis, Georgia Institute of
Technology.
\bibitem{phd175} B. B. Karklin$\ddot{\textrm{u}}$sh and V. B. Karkli$\check{\textrm{n}}$ (1976), {\it Inverse loops}, In
'Nets and Quasigroups', Mat. Issl. 39, 82–-86.
\bibitem{phd176} A. D. Keedwell and V. A. Shcherbacov (2002), {\it On m-inverse loops and
quasigroups with a long inverse cycle}, Australas. J. Combin. 26,
99–-119.
\bibitem{phd177} A. D. Keedwell and V. A. Shcherbacov (2003), {\it Construction and
properties of $(r, s, t)$-inverse quasigroups I}, Discrete Math.
266, 275–-291.
\bibitem{phd178} A. D. Keedwell and V. A. Shcherbacov, {\it Construction and properties
of $(r, s, t)$-inverse quasigroups II}, Discrete Math. 288 (2004),
61-–71.
\bibitem{phd33} M. K. Kinyon (2005), {\it A survey of Osborn loops},
Milehigh conference on loops, quasigroups and non-associative
systems, University of Denver, Denver, Colorado.
\bibitem{phd36} M. K. Kinyon, K. Kunen (2004), {\it The structure of
extra loops}, Quasigroups and Related Systems 12, 39--60.
\bibitem{phd35} M. K. Kinyon, K. Kunen, J. D. Phillips (2004), {\it
Diassociativity in conjugacy closed loops}, Comm. Alg. 32, 767--786.
\bibitem{phd89} J. M. Osborn (1961), {\it Loops with the weak
inverse property}, Pac. J. Math. 10, 295--304.
\bibitem{phd3} H. O. Pflugfelder (1990), {\it Quasigroups and loops : Introduction},
Sigma series in Pure Math. 7, Heldermann Verlag, Berlin, 147pp.
\bibitem{phd90} E. Wilson (1966), {\it A class of loops with the
isotopy-isomorphy property}, Canad. J. Math. 18, 589--592.
\end{thebibliography}
\end{document}